\documentclass[multi]{cambridge7A}
\usepackage[UKenglish]{babel}
\usepackage[longnamesfirst,sectionbib]{natbib}
\usepackage{chapterbib}
\usepackage{amsmath,amssymb,amsthm,amsfonts}
\usepackage{enumerate,url}

\setcitestyle{authoryear,round,semicolon}

\allowdisplaybreaks[1]

\providecommand*\Index[1]{#1\index{#1}}
\providecommand*\undex[1]{} 
\providecommand*\Undex[1]{#1} 

\begin{document}
\alphafootnotes
\author[Warren J. Ewens and Geoffrey A. Watterson]{Warren J.
  Ewens\footnotemark[1]$^\ast$
  and Geoffrey A. Watterson\footnotemark[2]}
\chapter[Kingman and mathematical population genetics]{Kingman and
  mathematical population genetics}
\footnotetext[1]{324 Leidy Laboratories, Department of Biology,
  University of Pennsylvania, Philadelphia, PA 19104, USA;
  wewens@sas.upenn.edu}
\footnotetext[2]{15 Brewer Road, Brighton East,
  Victoria 3187, Australia; geoffreywmailbox-monashfriends@yahoo.com.au}
\makeatletter
\def\@fnsymbol#1{*\relax }
\makeatother
\footnotetext[3]{Corresponding author}
\arabicfootnotes
\contributor{Warren J. Ewens
  \affiliation{University of Pennsylvania}}
\contributor{Geoffrey A. Watterson
  \affiliation{Monash University}}
\renewcommand\thesection{\arabic{section}}
\numberwithin{equation}{section}
\renewcommand\theequation{\thesection.\arabic{equation}}
\numberwithin{table}{section}
\renewcommand\thetable{\thesection.\arabic{table}}

\begin{abstract}
Mathematical population genetics is only one of Kingman's
many research interests. Nevertheless, his contribution to
this field has been crucial, and moved it in  several important
new directions. Here  we outline some aspects of his work
which have had a major influence on population genetics theory.
\end{abstract}

\subparagraph{AMS subject classification (MSC2010)}92D25

\section{Introduction}\index{mathematical genetics|(}
In the early years of the previous century, the main aim of population genetics
theory was to validate the Darwinian\index{Darwin, C. R.} theory
of \Index{evolution}, using the Mendelian\index{Mendel, G. J.}
hereditary\index{heredity} mechanism as the vehicle for determining how the characteristics
of any daughter generation depended on the corresponding characteristics of the
parental generation. By the 1960s, however, that aim had been achieved, and the theory
largely moved in a new, retrospective and statistical, direction.

This happened because, at that time,
\index{data!availability}data on the genetic constitution
of a population, or at least on a sample of individuals from that population, started to become available.
What could be inferred about the past history of the population  leading to these data? Retrospective questions
of this type include: ``How do we estimate the time at which \Index{mitochondrial Eve},
the woman whose
\index{DNA (deoxyribonucleic acid)!mitochondrial DNA}mitochondrial DNA is the most recent ancestor of the mitochondrial DNA currently
 carried in the human population, lived? How can contemporary genetic data be used
to track the `Out of Africa' migration? How do we detect signatures of past selective
events in our contemporary genomes\index{genome}?''
Kingman's\index{Kingman, J. F. C.!influence|(}  famous coalescent theory\index{Kingman, J. F. C.!Kingman coalescent}
became a central vehicle for addressing questions such as these. The very success of coalescent theory
has, however, tended to obscure Kingman's other contributions to
 population genetics theory. In
this note we review his various contributions to that theory, showing how coalescent theory
arose, perhaps naturally, from his earlier contributions.

\section{Background}
\label{back}
Kingman attended lectures in genetics at Cambridge in about 1960, and his earliest contributions
to population genetics date from 1961. It was well known at that time that in a randomly
mating population for which the fitness\index{fitness} of any individual depended on his genetic make-up at a single gene locus\index{gene!locus|(}, the mean fitness of the population increased from one generation to the next, or at least remained constant, if only two possible alleles\index{allele|(}, or gene types, often labelled $A_1$ and
$A_2$, were possible at that gene locus. However, it was well known that more than two alleles
could arise at some loci (witness the ABO \Index{blood group} system, admitting three possible alleles,
A, B and O).  Showing that in this case the mean population fitness is
non-decreasing in time under \Index{random mating} is far less easy to
prove. This was conjectured by Mandel and Hughes (1958)\index{Mandel,
S. P. H.}\index{Hughes, I. M.} and proved in the `symmetric'
case by Scheuer and Mandel (1959)\index{Scheuer, P. A. G.} and Mulholland and
Smith (1959)\index{Mulholland, H. P.}\index{Smith, C. A. B.},  and  more generally
 by Atkinson {\it et al.} (1960)\index{Atkinson, F. V.}\index{Moran, P. A. P.} and (very generally)
 Kingman, (1961a,b). Despite this success, Kingman then focused his research in areas quite different from genetics for the next fifteen years. The aim of this paper is to document some of his work following his re-emergence
into the genetics field, dating from 1976. Both of us were honoured to be associated with him in this work. Neither of us can remember the precise details, but the three-way interaction between the UK, the USA and Australia, carried out mainly
by the now out-of-date flimsy blue aerogrammes, must have started in  1976, and
continued during the time of Kingman's intense involvement in population genetics.
This note is a personal account, focusing on this interaction: many others were working in the field at the same time.

One of Kingman's research activities during the period 1961-1976 leads to our first `background' theme.
In 1974 he established (Kingman, 1975) a surprising and beautiful result,
found in the context of storage strategies\index{storage strategy}.  It is
well known that the symmetric $K$-dimensional Dirichlet
distribution\index{Dirichlet, J. P. G. L.!Dirichlet distribution}
\begin{equation}
\frac{\Gamma(K\alpha)}{\Gamma(\alpha)^K}(x_1x_2\cdots x_K)^{\alpha -1}\,
dx_1\,dx_2\ldots\, dx_{K-1},
\label{eqn:k1}
\end{equation}
where $x_i \geq 0$, $\sum x_j = 1$, does not have a non-trivial limit as $K \rightarrow \infty$, for
given fixed $\alpha$. Despite this, if we let $K \rightarrow \infty$ and $\alpha \rightarrow 0$
in such a way that the product $K\alpha$ remains fixed at a constant value
$\theta$, then the distribution of the \emph{order statistics}\index{order statistics} $x_{(1)}\geq x_{(2)} \geq x_{(3)} \geq \cdots$ converges to a non-degenerate limit. (The parameter $\theta$ will turn out to have an important
genetical interpretation, as discussed below.) Kingman called this the
Poisson--Dirichlet
distribution\index{Poisson, S. D.!PoissonDirichlet distribution@Poisson--Dirichlet distribution}, but we suggest that its true author be honoured and that it be
called the `Kingman distribution'\index{Kingman, J. F. C.!Kingman distribution}. We refer to
it by this name in this paper. So important has the distribution become in mathematics generally that a book has been written
devoted entirely to it
\index{Feng, S.}(Feng, 2010).
This distribution has a rather complex form, and aspects of this form are given below.

The Kingman distribution appears, at first sight, to have nothing to do with population genetics theory. However, as we show below, it turns out, serendipitously, to be central to that theory. To see why this is so, we turn to our second `background' theme, namely the development of
 population theory in the 1960s and 1970s.

The nature of the gene\index{gene|(} was discovered by Watson\index{Watson,
J. D.}
and Crick\index{Crick, F. H. C.} in 1953. For our purposes the most important
of their results is the fact that a gene is in effect a
\index{DNA (deoxyribonucleic acid)}DNA sequence of, typically, some 5000
bases, each \Index{base} being one of four types, A, G, C or T. Thus the
number of types, or alleles, of a gene consisting of 5000 bases is
$4^{5,000}$.  Given this number, we may for many practical purposes suppose
that there are infinitely many different alleles possible at any gene
locus\index{gene!locus|)}. However, gene sequencing\index{gene!sequencing}
methods took some time to develop, and little genetic information at the
fundamental DNA level was available for several decades after Watson and
Crick.

The first attempt at assessing the degree of genetic variation\index{genetic variation|(} from one person
to another in a population  at a less fundamental level depended on the technique  of \Index{gel electrophoresis},  developed in the 1960s. In loose terms, this method measures the electric charge on a gene, with the charge levels usually thought of as taking integer values only.
Genes having different electric charges
are of different allelic types, but it can well happen that genes of different allelic types have the same electric charge. Thus there is
 no one-to-one relation between charge level and allelic type.  A simple mutation\index{mutation|(} model assumes that a mutant
gene has a charge differing from that of its parent gene by either $\pm1$. We return to this model in a moment.

In 1974 Kingman travelled to Australia, and while there met Pat Moran (as it
happens, the  PhD supervisor of both authors of this paper), who was working
at that time on this `charge-state' model\index{mathematical genetics!charge-state model|(}. The two of them discussed
the properties of a stochastic  model involving a population of $N$ individuals, and hence $2N$ genes at any given locus.   The population is assumed to evolve by random \Index{sampling}:
any daughter generation of genes is found by sampling, with replacement, from the genes from the parent generation. (This is the well-known `Wright--Fisher' model of population genetics, introduced into the
population genetics literature independently by Wright (1931)\index{Wright,
S.} and Fisher (1922)\index{Fisher, R. A.}.)  Further, each daughter generation gene is assumed to inherit the same charge as that of its parent with probability $1-u$, and with probability $u$ is a charge-changing mutant, the change in charge being equally likely to be $+1$ and $-1$.

At first sight it might seem that, as time progresses, the charge
levels on the genes in future generations become dispersed over the entire array of positive and negative integers. But this is not so. Kingman recognized that there is a coherency to the
locations of the charges on the genes brought about by common ancestry and the genealogy
of the genes in any generation. In Kingman's words (Kingman 1976), amended here to our terminology, ``The probability that
 [two genes in generation $t$] have a common ancestor
gene\index{gene!ancestor gene} [in generation $s$, for $s < t$,] is
  $1 - (1-(2N)^{-1})^{t-s}$, which is near unity when $(t-s)$ is large compared to $2N$. Thus the
  [locations of the charges in any generation] form a coherent group, \ldots, and the relative distances between the [charges] remain \Index{stochastically bounded}''. We do not dwell here on the elegant theory that Kingman developed for this model,  and note only that
   in the above quotation we see here the beginnings of the idea of looking
   backward in time to discuss properties of genetic
   variation\index{genetic variation|)} observed in a contemporary
   generation. This viewpoint is central to Kingman's concept of the
   coalescent\index{Kingman, J. F. C.!Kingman coalescent}, discussed in detail below.

   Parenthetically, the question of the mean number of `alleles', or occupied
   charge states, in a population of size $N$ (2$N$ genes) is of some
   mathematical interest. This depends on the mutation rate\index{mutation|)} $u$ and the population size $N$. It was originally conjectured by Kimura
   and Ohta (1978)\index{Kimura, M.}\index{Ohta, T.} that this mean remains
   bounded as $N \rightarrow \infty$. However, Kesten (1980a,b)\index{Kesten, H.} showed that it
   increases indefinitely as  $N \rightarrow \infty$, but at an extraordinarily slow rate. More exactly, he found
   the following astounding result. Define $\gamma_0 = 1$, $\gamma_{k+1} =
   e^{\gamma_k}$, $k = 1$, 2, 3, \ldots,
   and $\lambda(2N)$ as the largest $k$ such that $\gamma_k < 2N$. Suppose that $4Nu= 0.2$. Then the random number of `alleles' in the
   population divided by $\lambda(2N)$ converges in probability to a constant whose value
   is approximately 2 as $N \rightarrow \infty$.  Some idea of the slowness of the divergence of the mean number of alleles can be found by observing that if
   $2N = 10^{1656520}$, then $\lambda(2N) = 3$.

  In a later paper (Kingman 1977a), Kingman extended the theory to the multi-dimensional case,
  where it is assumed that data are available on a vector of measurements on each gene. Much of
  the theory for the one-dimensional charge-state model carries through more or less immediately to the
  multi-dimensional case. As the number of dimensions increases, some of this theory established
  by Kingman bears on the `infinitely many alleles' model\index{mathematical genetics!infinitely-many-alleles model} discussed in the next paragraph, although as Kingman himself noted, the geometrical structure inherent
  in the model implies that a convergence of his results to those of the infinitely-many-alleles model
  does not occur, since the latter model has no geometrical structure.

  The infinitely-many-alleles model, introduced in the 1960s, forms the second  background development that we discuss. This model has two components. The first is a purely demographic, or genealogical, model of the population.
  There are many such models, and here we consider only the Wright--Fisher model referred to above.
  (In the contemporary literature many other such models are discussed in the
  context of the infinitely-many-alleles model, particularly those of Moran
  (1958)\index{Moran, P. A. P.} and Cannings (1974)\index{Cannings, C.},
  discussed in Section \ref{robust}.)
  The second component refers to the mutation\index{mutation|(} assumption, superimposed on this model. In the infinitely-many-alleles model this assumption is that any new mutant gene
  is of an allelic type never before seen in the population.  (This is motivated by the very large number of alleles  possible at any gene locus, referred to above.) The model also assumes
  that the probability that any gene is a mutant is some fixed value $u$, independent of the allelic
  type of the parent and of the type of the mutant gene.

  From a practical point of view, the model
 assumes a technology (relevant to the 1960s) which is able to assess whether any two genes are of the same or are of
  different allelic types (unlike the charge-state model\index{mathematical genetics!charge-state model|)}, which does not fully
  possess this capability), but which is not   able to distinguish any further
  between two genes
  (as would be possible, for example,  if
 the
\index{DNA (deoxyribonucleic acid)}DNA sequences of the
  two genes were known).
Further, since an entire generation of genes is
  never observed in practice, attention focuses on the allelic configuration
  of the genes in a sample of size $n$, where $n$ is assumed to be small compared to $2N$,
  the number of genes in the entire population.

  Given the nature of the mechanism assumed in this model for distinguishing the
  allelic types of the $n$ genes in the sample, the data in effect consist of
  a partition\index{random partition|(}
  of the integer $n$ described by the vector $(a_1, a_2, \ldots, a_n)$, where $a_i$ is the number of allelic types observed in the sample exactly $i$ times each. It is necessary that $\sum i a_i = n$, and it turns out that under this condition, and to a close approximation,  the stationary probability of observing this vector is
\begin{equation}
\frac{n!\theta^{\sum a_i}}{1^{a_1}2^{a_2}\cdots n^{a_n}a_1!a_2!\cdots a_n!S_n(\theta)} ,
\label{eqn:ke}
\end{equation}
where $\theta$ is defined as $4Nu$ and $S_n(\theta) = \theta(\theta+1)(\theta+2)\cdots (\theta+n-1),$
(Ewens (1972), Karlin and McGregor (1972))\index{Karlin, S.}\index{McGregor,
J. L.}.

The marginal distribution of the number $K = \sum a_i$ of distinct alleles
in the sample  is found from (\ref{eqn:ke}) as
\begin{equation}
\text{Prob}(K = k) = |S_n^k|\theta^k/S_n(\theta),
\label{eqn:kdist}
\end{equation}
where $S_N^k$ is a Stirling number\index{Stirling, J.!Stirling number} of
the first kind. It follows from (\ref{eqn:ke}) and (\ref{eqn:kdist}) that $K$ is a \Index{sufficient statistic} for $\theta$,
so that the conditional distribution of $(a_1, a_2, \ldots, a_n)$ given $K$ is independent of $\theta$.

 The relevance of this observation is as follows.  As noted
 above, the extent of \Index{genetic variation} in a population was, by
 electrophoresis\index{gel electrophoresis} and
 other methods, beginning to be understood in the 1960s. As a result
 of this knowledge,  and for reasons not discussed here, Kimura
 advanced (Kimura 1968)\index{Kimura, M.} the so-called `neutral
 theory'\index{mathematical genetics!neutral theory}, in which it was claimed
 that much of the genetic variation  observed
 did not have a selective basis. Rather, it was claimed that  it was the result
  of purely random changes in allelic frequency inherent in the random sampling evolutionary
  model\index{mathematical genetics!random-sampling model} outlined above. This (neutral) theory  then becomes
  the null hypothesis in a statistical testing procedure, with some selective
 mechanism being the alternative hypothesis. Thus the expression in
 (\ref{eqn:ke}) is the null hypothesis allelic-partition
 distribution\index{allele!allelic partition distribution} of the
  alleles in a sample of size $n$.  The fact that the conditional distribution of $(a_1, a_2, \ldots, a_n)$ given $K$
  is independent of $\theta$  implies that an objective  testing
  procedure for the neutral theory can be found free of unknown parameters.

  Both authors of this paper worked on aspects of this statistical testing theory during
  the period 1972--1978, and further reference to this is made below. The random sampling
  evolutionary scheme described above is no doubt a simplification of real
  evolutionary processes, so in order
  for the testing theory to be applicable to more general evolutionary models it is
natural to ask: ``To what extent does the expression in (\ref{eqn:ke}) apply for evolutionary models other than that described above?'' One of us (GAW) worked on this question in the mid-1970s (Watterson, 1974a, 1974b). This question is also discussed below.

\section{Putting it together}
\label{pit}
One of us (GAW) read Kingman's 1975 paper soon after it appeared and recognized
its potential application to population genetics theory. In the 1970s the joint density function (\ref{eqn:k1}) was well known to arise
in that theory when some fixed finite number $K$ of alleles is possible at
the gene locus of interest, with symmetric mutation\index{mutation|)} between these alleles.
In  population genetics theory one considers, as mentioned above, infinitely
many possible alleles at any gene locus, so that the relevance of Kingman's limiting
($K \rightarrow \infty$) procedure to the infinitely many alleles
model\index{mathematical genetics!infinitely-many-alleles model}, that
 is the relevance of the Kingman distribution\index{Kingman, J. F. C.!Kingman distribution},  became immediately
apparent.

This observation led (Watterson 1976) to a derivation of an explicit form for the joint
density function of the first $r$ order statistics\index{order statistics|(} $x_{(1)}$, $x_{(2)}$, \ldots , $x_{(r)}$ in
the Kingman distribution. (There is an obvious printer's error
in equation (8) of Watterson's paper.) This joint
density function was shown to be of the form
\begin{equation}
f(x_{(1)}, x_{(2)}, \ldots , x_{(r)}) = \theta^r \Gamma(\theta)e^{\gamma\theta}g(y)\{x_{(1)}x_{(2)}\cdots x_{(r)}\}^{-1}
x_{(r)}^{\theta -1},
\label{eqn:k2}
\end{equation}
where $y = (1-x_{(1)}- x_{(2)}- \cdots - x_{(r)})/x_{(r)}$, $\gamma$ is
Euler's constant\index{Euler, L.!Euler Mascheroni constant@Euler--Mascheroni constant}
$0.57721\ldots$,  and $g(y)$ is best defined
through the Laplace transform\index{Laplace, P.-S.!Laplace transform} equation
(Watterson and Guess (1977))\index{Guess, H. A.|(}
\begin{equation}
\int_0^{\infty} e^{-ty} g(y) dy = \exp\left(\theta \int_0^1 u^{-1} (e^{-tu}-1)\, du\right).
\label{eqn:k3}
\end{equation}
The expression (\ref{eqn:k2}) simplifies to
\begin{equation}
f(x_{(1)}, \ldots, x_{(r)}) = \theta^r \{x_{(1)}\cdots x_{(r)}\}^{-1} (1-x_{(1)}-\cdots -x_{(r)})^{\theta -1}
\label{eqn:k4}
\end{equation}
when $x_{(1)}+ x_{(2)}+ \cdots + x_{(r-1)} + 2x_{(r)} \geq 1$, and in particular,
\begin{equation}
f(x_{(1)}) = \theta (x_{(1)})^{-1} (1-x_{(1)})^{\theta -1}
\label{eqn:k5}
\end{equation}
when $\frac12 \leq x_{(1)} \leq 1$.

Population geneticists are interested in the probability
of `population monomorphism'\index{monomorphism}, defined in practice as the
probability that the most frequent allele arises in the population
with frequency in excess of 0.99.  Equation (\ref{eqn:k5})
implies that this probability is close to $1 - (0.01)^{\theta}$.

Kingman himself had placed some special emphasis on the largest of
the order statistics\index{order statistics|)}, which in the genetics context is the allele
frequency of the most frequent allele. This leads to interesting questions in
genetics. For instance,  Crow (1973)\index{Crow, J. F.} had asked: ``What is the probability that the most frequent allele in a population at any time is also the oldest allele in the
population at that time?'' A nice application of \Index{reversibility}
arguments for suitable population models allowed Watterson and Guess
(1977)\index{Guess, H. A.|)} to obtain a simple answer to this question.   In models where all alleles are equally fit, the probability that  any nominated allele will survive longest into the future is
(by a simple symmetry argument) its current frequency.
For time-reversible processes, this is also the probability that it is the oldest
allele in the population.  Thus conditional on the current allelic frequencies, the probability
that the most frequent allele is also the oldest is simply its frequency $x_{(1)}$. Thus the answer to Crow's question is simply the mean frequency of the most frequent allele. A formula for this mean frequency, as a function of the \Index{mutation} parameter $\theta$, together with some  numerical values, were  given in Watterson and Guess (1977), and a partial
listing is given in the first row of Table \ref{tmfao}. (We discuss the entries in the second row
of this table in Section \ref{GEM}.)

\begin{table}[h]\footnotesize
\begin{minipage}{106mm}
  \caption{Mean frequency of (a) the most frequent allele, (b) the oldest
  allele, in a population as a function of $\theta$. The probability that the
  most frequent allele is the oldest allele is also its mean frequency.}
  \label{tmfao}
\end{minipage}
\begin{tabular}{@{}ccccccccc@{}}\hline
$\theta$&0.1&0.2&0.5&1.0&2.0&5.0&10.0&20.0\\\hline
Most frequent&0.936&0.882&0.758&0.624&0.476&0.297&0.195&0.122\\
Oldest&0.909&0.833&0.667&0.500&0.333&0.167&0.091&0.048\\\hline
\end{tabular}
\end{table}

 As will be seen from the table,  the mean frequency $E(x_{(1)})$ of
 the most frequent allele decreases as $\theta$ increases.
Watterson and Guess (1977) provided the
bounds $(\frac12)^{\theta} \leq  E(x_{(1)}) \leq 1 - \theta(1-\theta) \log 2$, which give an idea of the value of $E(x_{(1)})$  for small values of $\theta$, and also showed that $E(x_{(1)})$ decreases asymptotically like (log $\theta$)/$\theta$, giving
 an idea of the value of $E(x_{(1)})$  for large $\theta$.

From the point of view of testing the neutral theory\index{mathematical
genetics!neutral theory} of Kimura\index{Kimura, M.}, Watterson (1977, 1978)
subsequently used properties of these  order statistics for testing the null
hypothesis that there are  no selective
forces\index{evolution!selective force of} determining observed allelic frequencies. He considered various alternatives, particularly heterozygote advantage\index{heterozygote advantage} or the presence of some deleterious alleles. For instance, in (Watterson 1977) he investigated the situation when all heterozygotes had a slight selective advantage over all homozygotes. The population truncated homozygosity $\sum_{i}^r x_i^2$  figures prominently in the allelic distribution corresponding to (\ref{eqn:k2}) and was thus studied as a test statistic for the null hypothesis of no selective advantage. Similarly, when only a random sample of $n$ genes is taken from the population, the sample homozygosity can be used as a test statistic
of neutrality.

Here we make a digression to discuss two of the values in the first row of Table \ref{tmfao}. It
is well known that in the case $\theta = 1$, the allelic partition formula (\ref{eqn:ke}) describes
the probabilistic structure of the lengths of the cycles in a \Index{random permutation} of the
numbers $\{1, 2, \ldots, n\}$. Each \Undex{cycle} corresponds to an allelic type and  in the
notation $a_j$ thus indicates the number of cycles of length $j$. Various  limiting ($n \rightarrow \infty$)
properties of random permutations have long been of interest (see for example
Finch (2003))\index{Finch, S. R.}. Finch (page 284) gives the limiting
mean of the normalized length of the longest cycle as $0.624\ldots$ in such a random
permutation, and this agrees with
the value listed in Table \ref{tmfao} for the case $\theta = 1$. (Finch also in effect gives the standard
deviation of this normalized length as $0.1921\ldots$.) Next, (\ref{eqn:k5}) shows
that the limiting probability that the (normalized) length of the longest cycle exceeds $\frac12$
is $\log 2.$ This is the limiting value of the exact probability   for  a random permutation of the
numbers $\{1, 2, \ldots, n\}$, which from (\ref{eqn:ke}) is $1 - \frac12 + \frac13 - \cdots \pm \frac1n$.

Finch also considers aspects of a \Index{random mapping} of $\{1, 2, \ldots, n\}$ to
$\{1, 2, \ldots, n\}$. Any such a mapping forms a random number of `components', each component
consisting of a \Undex{cycle} with a number (possibly zero) of branches
attached to it. Aldous (1985)\index{Aldous, D. J.} provides
a full description of these, with diagrams which help in understanding them.
Finch takes up the question of finding properties of the normalized size of the largest component of such a random mapping, giving (page 289) a limiting mean of $0.758\ldots$ for this. This agrees with the value in
Table \ref{tmfao} for the case $\theta = 0.5$. This is no coincidence: Aldous (1985) shows
that in a limiting sense (\ref{eqn:ke}) provides the limiting distribution of the
number and (unnormalized) sizes of the components of this mapping, with
now $a_j$  indicating the number of components of size $j$.  As a further result, (\ref{eqn:k5}) shows
that the limiting probability that the (normalized) size of the largest component of a random mapping  exceeds $\frac12$
is $\log (1+\sqrt{2}) \approx 0.881374$.

Arratia {\it et al.} (2003)\index{Arratia, R.}\index{Barbour, A. D.}\index{Tavar\'e, S.} show that
(\ref{eqn:ke}) provides, for various
values of $\theta$, the \Index{partition structure} of a variety of other
combinatorial objects for finite $n$, and presumably the Kingman
distribution\index{Kingman, J. F. C.!Kingman distribution|(} describes appropriate limiting $(n \rightarrow \infty)$
results. Thus the genetics-based equation  (\ref{eqn:ke}) and the Kingman distribution provide a unifying theme for these
objects.

The allelic partition  formula (\ref{eqn:ke}) was originally derived without
reference to the $K$-allele model\index{mathematical genetics!K allele model@$K$-allele model|(} (\ref{eqn:k1}), but was also found (Watterson, 1976)
from that model as follows.  We start with a population whose allele
frequencies are given by the Dirichlet distribution\index{Dirichlet, J. P. G. L.!Dirichlet distribution} (\ref{eqn:k1}). If a random sample of $n$ genes is taken from such a population, then given the population's allele frequencies, the sample allele frequencies have a multinomial distribution. Averaging this distribution over the population distribution (\ref{eqn:k1}), and then introducing the alternative order-statistic sample description $(a_1, a_2, \ldots, a_n)$ as above,  the limiting distribution is the partition formula (\ref{eqn:ke}), found by letting $K \rightarrow \infty$ and $\alpha \rightarrow 0$ in (\ref{eqn:k1}) in such a way that the product $K\alpha$ remains fixed at a constant value $\theta$.

\section{Robustness}
\label{robust}\index{robustness}
  As stated above, the expression (\ref{eqn:ke}) was first found by assuming a
  random sampling\index{sampling|(} evolutionary model\index{mathematical genetics!random-sampling model|(}.
  As also noted, it can also be arrived at by  assuming that a random sample
  of genes has been taken from an infinite population whose allele frequencies
  have the Dirichlet distribution (\ref{eqn:k1}). It applies, however, to
  further models. Moran (1958)\index{Moran, P. A. P.} introduced a `birth-and-death' model\index{mathematical genetics!birth-and-death model} in which, at each unit time point, a gene is chosen at random from the population to die.
   Another gene is chosen at random to reproduce. The new gene either inherits the allelic type of its parent (probability $1-u$), or is of a new allelic type, not so
   far seen in the population, with probability $u$. Trajstman
  (1974)\index{Trajstman, A. C.} showed that (\ref{eqn:ke}) applies as the stationary allelic partition distribution exactly for Moran's  model, but with $n$ replaced by the finite population number of genes $2N$ and with  $\theta$  defined as ${2Nu/(1-u)}$. More than this, if a random sample of size $n$ is taken without replacement from the Moran model population, it too has an exact description as in (\ref{eqn:ke}).  This result is a consequence of Kingman's (1978b) study of the consistency of the allelic properties of  sub-samples of samples. (In practice, of course, the difference between sampling with, or without, replacement is of little consequence for small samples from large populations.)
Kingman (1977a, 1977b) followed up this result by showing  that random sampling from various other population models,
including significant cases of the Cannings (1974) model\index{Cannings, C.!Cannings model}, could also be approximated by (\ref{eqn:ke}). This was important because several consequences of (\ref{eqn:ke}) could then be applied more generally than was first thought, especially for the purposes of testing of the neutral alleles postulate. He also used
the concept of
\index{noninterference@non-interference}`non-interference'
(see the concluding comments in Section \ref{partition}) as a further reason
for the robustness of (\ref{eqn:ke}).

\section{A convergence result}
\label{conv}
It was noted in Section \ref{pit} that Watterson (1976) was able to arrive at both
the Kingman distribution\index{Kingman, J. F. C.!Kingman distribution|)} and the
allelic partition formula (\ref{eqn:ke})\index{allele!allelic partition distribution} from
the same starting point (the `$K$-allele' model)\index{mathematical genetics!K allele model@$K$-allele model|)}. This makes it clear that there
must be a close connection between the two, and in this section we outline
Kingman's work (Kingman 1977b) which made this explicit. Kingman imagined a sequence of
populations in which the size of population $i$, ($i = 1$, 2, \ldots) tends to infinity as $i \rightarrow \infty$.
For any fixed $i$ and any fixed sample size $n$ of genes taken from the population, there will be some probability of the partition $\{a_1, a_2, \ldots, a_n\}$, where $a_j$ has the
definition given in Section \ref{back}.
Kingman then stated that this sequence of populations would have
the \emph{Ewens sampling property}\index{Ewens, W. J.!Ewens sampling property}
if, for each fixed $n$,  this corresponding sequence of probabilities
of $\{a_1, a_2, \ldots, a_n\}$ approached that given in (\ref{eqn:ke})
as $i \rightarrow \infty$. In a parallel fashion, for each fixed $i$ there will also be a probability
distribution for the \Index{order statistics} $(p_{1}, p_{2}, \ldots)$, where $p_{j}$ denotes
the frequency of the $j$th most frequent allele in the population. Kingman then stated that this sequence would have the \emph{Poisson--Dirichlet limit} if
this sequence of probabilities approached that given by the Poisson--Dirichlet
distribution\index{Poisson, S. D.!PoissonDirichlet distribution@Poisson--Dirichlet distribution}.
(We would replace `Poisson--Dirichlet'  in this sentence by `Kingman'.) He then showed
that this sequence of populations has the \emph{Ewens sampling property} if and only
if it has the Poisson--Dirichlet (Kingman distribution) limit.

The proof is quite technical and we do not discuss it here. We have noted that  the
Kingman distribution may be thought of as the distribution of the (ordered) allelic frequencies in an infinitely large population evolving as the random sampling infinitely-many-allele process,
so this result provides a beautiful (and useful) relation between population and
sample properties of such a population.

\section{Partition structures}
\label{partition}\index{partition structure|(}
By 1977 Kingman was in full flight in his investigation of various genetics problems. One line
of his work started with the probability distribution (\ref{eqn:ke}), and his initially innocent-seeming observation that the size $n$ of the sample of genes bears further
consideration. The size of a sample is generally taken in Statistics as being comparatively
uninteresting, but Kingman (1978b) noted that a sample of $n$ genes could be regarded as having arisen from a sample of $n+1$ genes, one of which was accidently lost, and that this observation
induces a consistency property on the probability of any partition of the number $n$. Specifically,
he observed that if we write $\text{P}_n(a_1, a_2, \ldots )$ for the probability of the sample partition in a sample of size $n$, we require
\begin{align}
\text{P}_n(a_1, a_2, \ldots )
 &= \frac{a_1 + 1}{n+1}\text{P}_{n+1}(a_1 +1, a_2, \ldots ) +{}\notag \\
 &\quad\sum_{j=2}^{n+1}\frac{j(a_j + 1)}{n+1}\text{P}_{n+1}(a_1, \ldots,a_{j-1}-1,a_j +1,\ldots ).
\label{pstructure}
\end{align}
Fortunately, the distribution (\ref{eqn:ke}) does satisfy this equation. But
Kingman went on to ask a deeper question: ``What are the most general distributions that satisfy  equation (\ref{pstructure})?'' These distributions
he called `partition structures'.  He showed  that all such distributions
that are of interest in genetics could be represented in the form
\begin{equation}
\text{P}_n(a_1, a_2, \ldots ) = \int \text{P}_n(a_1, a_2, \ldots
  |\textbf{x})\,\mu(d\textbf{x})
\label{repmeasure}
\end{equation}
where $\mu$ is some probability measure over the space of infinite sequences $(x_1, x_2, x_3 \ldots)$ satisfying
$x_1 \geq x_2 \geq x_3 \cdots$,  $\sum_{n=1}^{\infty} x_n = 1$.

An intuitive understanding of this equation is the following. One way to obtain a consistent set of
distributions satisfying (\ref{pstructure}) is to imagine a hypothetically infinite
population of types, with a proportion $x_1$ of the most frequent type, a proportion $x_2$ of
the second most frequent type, and so on, forming a vector {\bf{x}}. For a fixed value of $n$, one could then imagine
taking a sample of size $n$ from this population, and write $\text{P}_n(a_1, a_2, \ldots |\textbf{ x})$ for
the (effectively multinomial) probability that the configuration of the sample is $(a_1, a_2, \ldots)$. It
is clear that the resulting sampling probabilities will automatically satisfy the consistency property
in (\ref{pstructure}). More generally one could imagine the composition of the infinite population
itself being random, so that first one chooses its composition {\bf{x}} from  $\mu$,
and then conditional on {\bf{x}} one takes a sample of size $n$ with probability $\text{P}_n(a_1, a_2, \ldots |\textbf{x})$. The right-hand side in (\ref{repmeasure}) is then the probability of obtaining the sample
configuration $(a_1, a_2, \ldots)$ averaged over the composition of the population. Kingman's remarkable result was
that all partition structures arising in genetics must have the form (\ref{repmeasure}), for some $\mu$.
Kingman called partition structures that could be
expressed as in (\ref{repmeasure}) `representable partition structures' and $\mu$ the `representing
measure', and later (Kingman
1978c)  found a representation generalizing (\ref{repmeasure}) applying
for any partition structure.

The similarity between (\ref{repmeasure}) and the celebrated de
Finetti\index{Finetti, B. de} representation
theorem for exchangeable\index{exchangeability} sequences might be noted. This has been explored by
Aldous (1985)\index{Aldous, D. J.|(} and Kingman (1978a), but we do not pursue the details
of this here.

In the genetics context,
the results of Section \ref{robust} show that samples from
Moran's\index{Moran, P. A. P.} infinitely many neutral alleles model\index{mathematical genetics!infinitely-many-alleles model}, as well as the population as a whole, have the partition structure
property. So do samples of genes from other genetical models. This makes it
natural to ask: ``What is  the representing measure $\mu$ for the allelic
partition distribution\index{allele!allelic partition distribution} (\ref{eqn:ke})?'' And
here we come full circle, since he showed that the required representing
measure is the Kingman distribution\index{Kingman, J. F. C.!Kingman distribution|(}, found by him in (Kingman, 1975) in quite a different context!

 The relation between the Kingman distribution and the
  sampling distribution (\ref{eqn:ke}) is of course connected to the convergence results
discussed in the previous section.  From the point of view of the geneticist,
the Kingman distribution is then regarded as applying for an infinitely large population,
evolving essentially via the random sampling\index{sampling|)} process\index{mathematical genetics!random-sampling model|)} that led to (\ref{eqn:ke}). This was
made precise by Kingman in (1978b), and it makes
it unfortunate that the Kingman distribution does not have a `nice' mathematical form.
However, we see in Section \ref{GEM} that a very pretty analogue of the Kingman distribution exists when we label alleles not by their frequencies but by their ages in the population.
This in turn  leads to the capstone of Kingman's work in genetics, namely the
  coalescent process\index{Kingman, J. F. C.!Kingman coalescent}.

Before discussing these matters we mention another  property enjoyed by the distribution
(\ref{eqn:ke}) that Kingman investigated, namely that of
\index{noninterference@non-interference}non-interference. Suppose
that we take a gene at random from the sample of $n$ genes, and find that there are in all
$r$ genes of the allelic type of this gene in the sample. These $r$ genes are
now removed, leaving $n-r$ genes. The non-interference requirement is that
the probability structure of these $n-r$ genes should be the same as
that of an original sample of $n-r$ genes, simply replacing $n$ wherever
found by $n-r$. Kingman showed that of all partition structures of interest
in genetics, the only one also satisfying this non-interference requirement is
(\ref{eqn:ke}).  This explains in part the \Index{robustness} properties of (\ref{eqn:ke}) to
various evolutionary genetic models. However, it also has a natural
interpretation in terms of the coalescent process, to be discussed in Section \ref{coal}.

We remark in conclusion that the partition structure concept has become
influential not only in the genetics context, but in Bayesian
statistics\index{Bayes, T.!Bayesian statistics},
mathematics and various areas of science, as the papers of Aldous
(2009)\index{Aldous, D. J.|)} and
of Gnedin, Haulk and Pitman (2009)\index{Gnedin, A. V.}\index{Haulk, C.}\index{Pitman, J. [Pitman, J. W.]} in this Festschrift show.  That this should
be so is easily understood when one considers the natural logic of the ideas
leading to it.\index{random partition|)}\index{partition structure|)}

\section{`Age' properties and the GEM distribution}
\label{GEM}
We have noted above that the Kingman distribution is not user-friendly.
This makes it all the more
interesting that a \emph{size-biased}\index{size biasing} distribution closely related to it,
namely the GEM
distribution\index{GEM distribution|(}, named for Griffiths (1980), Engen (1975) and McCloskey (1965), who
established its salient properties,
is both simple and elegant, thus justifying the acronym `GEM'. More important, it has a central interpretation
with respect to the ages of the alleles in a population. We now describe this distribution.

We have shown that the ordered allelic frequencies in the population follow the
Kingman  distribution.
Suppose that a gene is taken at random from the population.  The probability that this gene will be
of an allelic type whose frequency in the population is $x$ is just $x$.
This allelic type was thus sampled by this choice in a
size-biased\index{sampling!size-biased sampling}
way. It can be shown from properties of the Kingman
distribution\index{Kingman, J. F. C.!Kingman distribution|)}
that the probability density  of the frequency of the allele
determined by this randomly chosen gene is
\begin{equation}
f(x) = \theta (1-x)^{\theta - 1}, \quad 0  < x < 1.
\label{eqsbsbyp}
\end{equation}
This result was also established by Ewens (1972).

Suppose now that all genes of the allelic type just chosen are removed
from the population. A second gene
is now drawn at random from the population and its allelic type observed. The frequency of the allelic type of this gene
among the genes remaining at this stage is also given by (\ref{eqsbsbyp}).
All genes of this second allelic type  are now also removed from the population.
A third gene then drawn at random from the genes remaining, its allelic type observed, and all genes
of this (third) allelic type removed from the population. This process is continued
indefinitely.
At any stage, the distribution
of the frequency of the allelic type of any gene just drawn among
the genes left when the draw takes place is given by (\ref{eqsbsbyp}).
This leads to the following representation. Denote by $w_j$ the population
frequency of the $j$th allelic type drawn. Then we can write
\begin{equation}
w_1 = x_1,\  \ldots,\ w_j = (1-x_1)(1-x_2)\cdots (1-x_{j-1})x_j,\quad (j = 2, 3, \ldots ),
\label{eqformforgem}
\end{equation}
where the $x_j$ are independent random variables, each having the
distribution (\ref{eqsbsbyp}). The random vector $(w_1, w_2, \ldots )$ then
has the GEM distribution.

All the alleles in the population at any time eventually leave the population,
through the joint processes of \Index{mutation} and random drift, and any allele with
current population frequency $x$ survives the longest with probability $x$.
That is, since the GEM distribution was found according to a size-biased process, it
also arises when alleles are labelled
according to the length of their
future persistence in the population.
Time reversibility arguments then show that the GEM distribution also applies when
the alleles in the population are labelled by their age. In other words, the
vector $(w_1, w_2, \ldots )$ can be thought of as the vector of allelic
frequencies when alleles are ordered with respect to their ages in the population
(with allele 1 being the oldest).

The Kingman coalescent\index{Kingman, J. F. C.!Kingman coalescent}, to be discussed in the following section, is concerned among
other things with `age' properties of the alleles in the population. We
thus present some of these properties here as an introduction
to the coalescent: a more complete list can be found in Ewens (2004).
The elegance of many age-ordered formulae derives
directly from the simplicity and tractability  of the GEM distribution.

Given the focus on retrospective questions, it is
natural to ask questions about the oldest allele in the population.
The GEM distribution shows that the mean population frequency
 of the oldest allele in the population is
\begin{equation}
\theta \: \int_0^1 x (1-x)^{\theta - 1} \,dx = \frac1{1+\theta}.
\label{issswrp}
\end{equation}
This implies that when $\theta$ is very small, this mean frequency is approximately
$1 - \theta$. It is interesting to compare this with the
mean frequency of the most frequent allele when $\theta$ is
small, found in effect from the Kingman distribution\index{Kingman, J. F. C.!Kingman distribution} to be
approximately $1 - \theta \log 2$.  A more general
set of comparisons of these two mean frequencies, for
representative values of $\theta$,
is given in Table \ref{tmfao}.

More generally,
the mean population frequency of the $j$th oldest allele in the population is
$$\frac{1}{1+\theta} \Bigl(\frac{\theta}{1+\theta}\Bigr)^{j-1}.$$
For the case $\theta = 1$, Finch (2003)\index{Finch, S. R.} gives the mean frequencies
of the second and third most frequent alleles as $0.20958\ldots$
and $0.088316\ldots$ respectively, which may
 be compared to the mean frequencies of the second and third oldest alleles,
 namely $0.25$ and $0.125$.  For $\theta = 1/2$
the mean frequency of the second most
frequent allele is  $0.170910\ldots$, while the mean frequency
of the second oldest allele is $0.22222$.

Next, the probability that a gene drawn at random from the population is
of the type of the oldest allele is the mean frequency of the oldest allele, namely
$1/(1+\theta),$ as just shown (see also Table \ref{tmfao}). More generally the probability that $n$ genes drawn at random
from the population are all of the type of the oldest allele  in the population is
\begin{equation}
 \theta \int_0^1 x^{n} (1-x)^{\theta - 1} \, dx  = \frac{n!}{(1+\theta) (2+\theta)\cdots  (n+\theta)}.
 \label{eqn.freqoldest}
 \end{equation}

The GEM distribution  has a number
of interesting mathematical properties, of which we mention
here only one.
It is a so-called `residual allocation' model\index{residual allocation model} (Halmos 1944)\index{Halmos, P. R.}.
Halmos envisaged a king with one kilogram of gold dust, and
an infinitely long line of beggars asking for gold. To the first
beggar the king gives $w_1$ kilogram of gold, to
the second  $w_2$ kilogram of gold, and so on,
as specified in (\ref{eqformforgem}),
where the $x_j$ are independently and identically distributed (i.i.d.)\ random variables, each having
some probability distribution over the interval $(0,1)$.

Different forms
of this distribution lead to different properties of the distribution of
the  `residual allocations' $w_1$, $w_2$, $w_3$, \ldots.
One such property is that the distribution of $w_1$, $w_2$, $w_3$, \ldots\
be invariant under size-biased sampling\index{sampling!size-biased sampling}. It can be shown that the GEM
distribution is the only residual allocation model having this property. This
fact had been exploited by
Hoppe (1986, 1987)\index{Hoppe, F.}
to derive various results of interest in genetics and ecology.

We now turn to sampling\index{sampling|(} results. The probability that $n$ genes drawn at random from the population   are all of the same
allelic type as the oldest allele in the population is given in (\ref{eqn.freqoldest}).
The probability that $n$ genes drawn at random from the population   are all of the same
unspecified allelic type is
 $$\theta \int_0^1 x^{n-1} (1-x)^{\theta - 1} \, dx  = \frac{(n-1)!}{(1+\theta) (2+\theta)\cdots  (n+\theta -1)},$$
in agreement with (\ref{eqn:ke}) for the case $a_j = 0$, $j = 1$, 2, \ldots, $n-1$, $a_n = n$.
From this result and that in (\ref{eqn.freqoldest}), given that $n$ genes drawn at random are all of the same allelic type, the probability that
they are all of the allelic
type of the oldest allele   is $n/(n +\theta)$. The similarity of this expression with
that deriving from a Bayesian
calculation\index{Bayes, T.!Bayesian statistics} is of some interest.

 Perhaps the most important sample distribution concerns the frequencies of
 the alleles in the sample when ordered by age. This
 distribution was found by
Donnelly and  Tavar\'{e} (1986)\index{Donnelly, P. [Donnelly, P. J.]}\index{Tavar\'e, S.},
who showed that the probability that the number of alleles in the sample takes
the value $k,$ and that the age-ordered numbers
of these alleles in the sample are, in age order, $n_{(1)}$,  $n_{(2)}$, \ldots , $n_{(k)}$, is
\begin{equation}
\frac{\theta^k (n-1)!}{S_n(\theta) n_{(k)}(n_{(k)}+ n_{(k-1)}) \cdots (n_{(k)}+ n_{(k-1)} + \cdots n_{(2)})},
\label{dontav}
\end{equation}
where $S_j(\theta)$ is defined below (\ref{eqn:ke}). This formula can be found
in several ways, one being as the size-biased\index{size biasing} version of  (\ref{eqn:ke}).

These are many interesting  results connecting the oldest allele in the sample to the oldest
allele in the population.   For example, Kelly (1976)\index{Kelly, F. P.}
showed that the probability that the oldest allele in the sample is represented
$j$ times in the sample is
\begin{equation}
\frac{\theta}{n} \binom{n}{j} \binom{n+\theta - 1}{j}^{-1}, \quad j =  1, 2, \ldots, n.
\label{petergorold}
\end{equation}
He also showed that the probability that the oldest allele in the population is observed at all
in the sample
is $n/(n+\theta).$ The probability that a gene seen $j$ times
in the sample is of the oldest   allelic type in the population is
  $j/(n+\theta).$   When $j=n$, so that there is only one
  allelic type present in the sample, this probability is $n/(n + \theta)$.
Donnelly  (1986) showed, more generally, that the probability that the
oldest allele in the population is observed $j$
times in the sample is
\begin{equation}
\frac{\theta}{n + \theta} \binom{n}{j} \binom{n+\theta - 1}{j}^{-1}, \quad j = 0, 1, 2, \ldots, n.
\label{unconprobo}
\end{equation}
This is of course closely connected to
Kelly's result.  For the case $j = 0$   the probability (\ref{unconprobo})   is $\theta/(n+\theta)$, confirming
the complementary probability  $n/(n+\theta)$ found above.     Conditional on the event
that the oldest allele in the population does appear in the sample, a straightforward
calculation using (\ref{unconprobo})    shows that  this   conditional
probability and that in   (\ref{petergorold}) are identical.

It will be expected that various
exact results hold  for the Moran\index{Moran, P. A. P.} model, with $\theta$ defined as $2Nu/(1-u)$.
The first of these is an exact representation
of the GEM distribution, analogous to (\ref{eqformforgem}). This has been
provided by Hoppe (1987)\index{Hoppe, F.}. Denote by $N_1$, $N_2$, \ldots\ the numbers of genes of the
oldest, second-oldest, $\ldots$ alleles in the population.
Then $N_1$, $N_2$, \ldots\ can be defined in turn by
\begin{equation}
N_i = 1 + M_i,    \quad i = 1, 2, \ldots,
\label{eqbbyypls}
\end{equation}
where $M_i$ has a binomial distribution with index $2N - N_1 - N_2 - \cdots - N_{i-1} -1$
and parameter $x_i$, where $x_1$, $x_2$, \dots\ are i.i.d.\ continuous random variables each
having the density  function (\ref{eqsbsbyp}).
      Eventually $N_1 + N_2 +\cdots+N_k = 2N$
and the process stops, the final index $k$ being identical  to the number $K_{2N}$ of alleles
in the population.

It follows directly from this representation that the mean of $N_1$ is
$$1 + (2N-1)\theta \int_0^1 x (1-x)^{\theta - 1} \, dx \; = \; \frac{2N+\theta}{1+\theta}.$$

If there is only one allele in the population, so that the population is strictly
monomorphic\index{monomorphism}, this allele
must be the oldest one in the population.     The above representation shows that the probability that the oldest
allele arises $2N$ times in the population is
$$\text{Prob\/} \, (M_1 = 2N-1) = \theta \int_0^1 x^{2N-1} (1-x)^{\theta - 1} \, dx,$$
and this reduces to the exact monomorphism probability
$$ \frac{2N-1}{(1+\theta)(2+\theta)\cdots(2N-1+\theta)}$$
for the Moran model.

More generally, Kelly (1977)\index{Kelly, F. P.}
has shown that the probability that  the oldest allele in the population
is represented by $j$ genes is, exactly,
\begin{equation}
\frac{\theta}{2N}\binom{2N}{j}\binom{2N +\theta -1}{j}^{-1}.
\label{kellygorold}
\end{equation}
The case $j = 2N$ considered above is a particular example of (\ref{kellygorold}), and
the mean number $(2N + \theta)/(1+\theta)$ also follows from (\ref{kellygorold}).

We now consider `age' questions.  It is found that
 the mean time, into the past,
that the oldest allele in the population entered the population (by a \Index{mutation} event) is
\begin{equation}
\label{eq8.2asold}
\text{Mean age of oldest allele\/}\; = \; \sum^{2N}_{j =1} \frac{4N}{j(j +\theta -1)} \;\text{generations}.
\end{equation}
It can  be shown  (see Watterson and Guess (1977)\index{Guess, H. A.} and
Kelly (1977)\index{Kelly, F. P.})
that not only the mean age of the oldest allele, but indeed the entire probability distribution of
its age, is independent of its current frequency and indeed of the
 frequency of all alleles in the population.

If an allele is observed in the
population with frequency $p,$ its mean age is
\begin{equation}
\sum_{j=1}^{2N} \frac{4N}{j(j+\theta -1)}\Bigl(1 - (1-p)^j\Bigr) \; \text{generations}.
\label{eq3.20newnew}
\end{equation}
This is a generalization of the expression in (\ref{eq8.2asold}), since if $p=1$
only one allele exists in the population, and it must then be the oldest allele.

Our final calculation concerns the mean age of the oldest allele in a sample
of $n$ genes. This is
\begin{equation}
4N \sum_{j = 1}^n \frac{1}{j(j+\theta - 1)} \;\text{generations}.
\label{meansampleage}
\end{equation}
Except for small values of $n$, this is close to the mean age of the oldest
allele in the population, given in (\ref{eq8.2asold}). In other words, unless $n$ is small,
it is likely that the oldest allele in the population is represented in the sample.

We have listed the various results given in this section not only because of their
intrinsic interest, but because they form a natural lead-in to Kingman's celebrated
coalescent process\index{Kingman, J. F. C.!Kingman coalescent|(}, to which we now turn.\index{GEM distribution|)}

\section{The coalescent}
\label{coal}
The concept of the coalescent is now discussed at length in many textbooks, and entire books (for example
Hein, Schierup and Wiuf (2005)\index{Hein, J.}\index{Schierup,
M. H.}\index{Wiuf, C.} and
Wakeley (2009))\index{Wakeley, J.} and book chapters (for example Marjoram and
Joyce (2009)\index{Marjoram, P.}\index{Joyce, P.} and Nordborg
(2001))\index{Nordborg, M.}
have been written about it.  Here we can do no more than outline the salient aspects
of the process.

The aim of the coalescent is to describe the common
ancestry\index{ancestral history} of the sample of $n$ genes at various times in the past through the concept of
an equivalence class\index{equivalence class/relation|(}. To do this we introduce the notation $\tau$, indicating a time $\tau$
in the past (so that if $\tau_1 > \tau_2$, time $\tau_1$ is further in
the past than time $\tau_2$). The sample of $n$ genes is assumed taken
at time $\tau = 0$.

  Two genes in the sample of $n$ are in the same equivalence class  at time $\tau$
if they have a common ancestor at this time. Equivalence classes
are denoted by parentheses: Thus if $n = 8$ and at time $\tau$ genes 1 and 2 have one common
ancestor,   genes 4 and 5 a second, and genes   6 and 7 a third, and
none of the three common ancestors are identical and none is
identical to the ancestor of gene 3 or of gene 8 at time $\tau$,  the equivalence
classes at time $\tau$ are
\begin{equation}
\{(1,2), (3), (4,5), (6,7), (8)\}.
\label{equivclases}
\end{equation}
We call any such set of equivalence classes an equivalence relation, and
denote any such equivalence relation by a Greek letter. As two particular
cases, at time $\tau = 0$ the equivalence relation is $\phi_1 =
\{(1),(2),(3),(4),(5),\break (6),(7),(8)\}$, and at the time of the most recent
common ancestor of all eight genes, the equivalence relation is $\phi_n =
\{(1,2,3,4,5,6,7,8)\}$.  The Kingman coalescent process is a description of
the details of the ancestry of the $n$ genes moving from $\phi_1$ to $\phi_n$.
For example, given the equivalence relation in (\ref{equivclases}), one
possibility for the equivalence relation following a coalescence is $\{(1,2),
(3), (4,5), (6,7,8)\}$.  Such an amalgamation is called a coalescence, and the
process of successive such amalgamations is called the coalescence process.

Coalescences are assumed to take place according to a Poisson
process\index{Poisson, S. D.!Poisson process}, but
with a rate depending on the number of equivalence classes present. Suppose that there are $j$ equivalence classes at time $\tau$. It is assumed that
no coalescence takes places between time $\tau$ and  time $\tau + \delta \tau$ with probability
 $1 - \frac12 j(j-1) \delta \tau$. (Here and throughout we ignore terms
of order $(\delta \tau)^2$.)  The probability
that the process moves from one nominated equivalence class (at time $\tau$)
to some nominated equivalence class which can be derived from it is
 $\delta \tau$. In other words, a coalescence takes place in this time interval
 with probability  $ \frac12 j(j-1) \delta \tau$, and all of the $j(j-1)/2$  amalgamations
 possible at time $\tau$ are equally likely to occur.

In order for this  process to describe the `random sampling' evolutionary model\index{mathematical genetics!random-sampling model|(} described
above, it is necessary to scale time so that unit time corresponds to $2N$
generations. With this scaling, the time $T_j$ between  the formation of  an
equivalence relation with $j$ equivalence classes to one with $j-1$ equivalence
classes has an exponential distribution with mean $2/j(j-1)$.

The (random) time $T_{\text {MRCAS}} = T_n + T_{n-1} +
T_{n-2} + \cdots + T_2$ until all genes in  the sample first had just one common ancestor    has mean
\begin{equation}
E(T_{\text {MRCAS}}) = 2\sum_{j=2}^n \frac{1}{j(j-1)} =
 2\left(1-\frac{1}{n}\right).
\end{equation}
(The suffix `MRCAS' stands for `most recent common ancestor of
the sample.)
This is, of course close to 2 coalescent time units, or 4N generations, when
$n$ is large. Tavar\'{e} (2004)\index{Tavar\'e, S.} has found the (complicated)
distribution of $T_{\text {MRCAS}}$.  Kingman (1982a,b,c) showed that for large populations, many
population models (including the `random sampling' model)\index{mathematical genetics!random-sampling model|)} are well approximated in their sampling\index{sampling|)}
attributes by the coalescent process. The larger the population the more accurate is this approximation.

We now introduce mutation\index{mutation|(} into the coalescent. Suppose that the probability that any
particular ancestral gene mutates in the time interval $(\tau+\delta\tau, \tau)$ is
$\frac{\theta}{2}\delta\tau$. All mutants are assumed to be of new allelic types (the
infinitely many alleles assumption).
If at time $\tau$ in the coalescent there are $j$ equivalence
classes\index{equivalence class/relation|)}, the probability
that either a mutation or a coalescent event had occurred in $(\tau+\delta\tau, \tau)$ is
\begin{equation}
 j\frac{\theta}{2}\delta\tau + \frac{j(j-1)}{2} \delta\tau = \frac{1}{2}j(j+\theta-1)\delta\tau .
\end{equation}
We call such an occurrence a defining event, and given that a defining event did occur, the probability
that it was a mutation is $\theta / (j + \theta -1)$ and that it is a coalescence is $(j-1) /(j + \theta -1)$.

The probability that $k$ different allelic types are seen in the sample is then the
probability that $k$ of these defining events were mutations. The above reasoning shows that
this probability must be proportional to $\theta^k/S_n(\theta)$, where $S_n(\theta)$ is defined below
(\ref{eqn:ke}), the constant of proportionality being independent of $\theta$. This argument leads to
(\ref{eqn:kdist}).

Using these results and  combinatorial arguments counting all possible coalescent paths
from a partition $(a_1, a_2, \ldots, a_n)$ back
to the original common ancestor, Kingman (1982a) was able to derive the more detailed sample partition probability
distribution (\ref{eqn:ke}), and deriving this distribution from coalescent arguments
is perhaps  the most pleasing way of arriving at it. For further comments along these lines, see (Kingman (2000)).

The description of the coalescent given above follows the original derivation given by  Kingman (1982a).
The coalescent is perhaps more naturally understood as a random binary tree. These have now been investigated in great detail:
see for example Aldous and Pitman (1999)\index{Aldous, D. J.}\index{Pitman, J. [Pitman, J. W.]}.

Many genetic results can be obtained quite simply by using the coalescent ideas. For example,
Watterson and Donnelly (1992)\index{Donnelly, P. [Donnelly, P. J.]} used Kingman's coalescent to discuss the
question ``Do Eve's Alleles Live On?''\index{mitochondrial Eve}
To answer this question we assume the infinitely-many-neutral-alleles
model\index{mathematical genetics!neutral theory} for the population and consider a random sample of $n$ genes taken at time `now'.
 Looking back in time, the ancestral lines of those genes coalesce to the MRCAS, which may be
called the sample's `Eve'. Of course if Eve's allelic type survives into the sample it would be the oldest,
but it may not have survived because of intervening mutation\index{mutation|)}. If we denote by $X_{n}$ the number of
representative genes of the oldest allele, and by $Y_{n}$ the number of genes having Eve's allele, then
Kelly's\index{Kelly, F. P.} result (\ref{petergorold}) gives the distribution of $X_{n}$.  We denote that distribution
here by $p_{n}(j)$, $j = 0$, 1, 2, \ldots, $n$, and the distribution of $Y_{n}$ by
$q_{n}(j)$, $j = 0$, 1, 2, \ldots, $n$. Unlike the simple explicit expression for $p_{n}(j)$, the
corresponding expression for $q_{n}(j)$ is very complicated: see (2.14) and (2.15) in Watterson
and Donnelly (1992), derived using some of Kingman's (1982a) results. Using the relative probabilities
of a mutation or a coalescence at a defining event gives rise to a recurrence equation for
$q_{n}(j)$, $j = 0$, 1, 2, \ldots, $n$ as
\begin{align}
&[n(n-1) +j\theta]q_{n}(j)\notag \\
 &\quad{}= n(j-1)q_{n-1}(j-1)
+n(n-j-1) q_{n-1}(j) +(j+1)\theta q_{n}(j+1)
\label{Griffiths}
\end{align}
for $j = 0$, 1, 2, \ldots, $n$, (provided that we interpret $q_{n}(j)$
as zero outside this range),  and for $n = 2$, 3, \ldots.  The boundary conditions $ q_{1}(j) =1$ for
$ j = 1$ , $q_1(j) = 0 $ for $ j >  1$, and
$$q_{n}(n) = p_{n}(n) = \prod_{k=2}^n\frac{k-1}{k+ \theta - 1}$$
apply, the latter  because if $X_{n} = n$ then all sample genes descend from a gene having the oldest allele,
and `she' must be Eve. The recurrence (\ref{Griffiths}) is a special case of one found
by Griffiths (1989)\index{Griffiths, R. C.} in his equation (3.7).

The expected number of genes of Eve's allelic type was given by Griffiths
(1986), (see also Beder (1988))\index{Beder, B.}, as
\begin{equation}
\text{E}(Y_{n}) = \sum_{j=0}^n jq_{n}(j) = n\prod_{j=2}^n\frac{j(j-1)}{j(j-1)+\theta}.
\end{equation}
Watterson and Donnelly (1992)\index{Donnelly, P. [Donnelly, P. J.]} gave some numerical examples, some asymptotic results, and
some bounds for the distribution $ q_{n}(j)$, $j = 0$, 1, 2, \ldots, $n$. One result of interest is
that $q_{n}(0)$, the probability of Eve's allele being extinct in the sample, increases with $n$,
to $q_{\infty}(0)$ say. One reason for this is that a larger sample may well have its `Eve' further back in the past than a
smaller sample. We might interpret $q_{\infty}(0)$ as being the probability that an infinitely large population has lost its `Eve's' allele\index{allele|)}. Note that the bounds
\begin{equation}
\frac{{\theta}^2}{(2+\theta)(1+\theta)} < q_{\infty}(0) \le \frac{\theta e^{\theta} - \theta}{\theta e^{\theta} +1},
\end{equation}
for $ 0 < \theta < \infty$, indicate that for all  $\theta$ in this range, $ q_{\infty}(0) $ is neither 0 nor 1. Thus, in contrast to
the situation in  branching processes\index{branching process}, there are no sub-critical or super-critical phenomena here.\index{Kingman, J. F. C.!Kingman coalescent|)}

\section{Other matters}
\label{other}
There are many other topics that we could mention in addition to those described above.
On the mathematical side, the Kingman distribution has a close connection to
prime factorization\index{prime factorization} of large integers. On the genetical side, we have not mentioned
the `infinitely many sites'
model\index{mathematical genetics!infinitely-many-sites model}, now
frequently used by geneticists, in which the
\index{DNA (deoxyribonucleic acid)}DNA structure of the gene\index{gene|)}
plays a central role. It is a tribute to Kingman that his
work opened up more topics than can be discussed here.

\paragraph{Acknowledgements}
Our main acknowledgement is to John Kingman himself.
The power and beauty of his work was, and still is, an inspiration to us both.
His generosity, often ascribing to us ideas of his own,
was unbounded. For both of us,
working with him was an experience never to be forgotten. More generally the
field of population genetics owes him an immense and, fortunately, well-recognized  debt.
We also thank an anonymous referee for suggestions which substantially improved this paper.\index{mathematical genetics|)}\index{Kingman, J. F. C.!influence|)}

\end{document}